\documentclass[oneside,12pt]{amsart}
\usepackage{amssymb, amscd}
\usepackage[all]{xy}
\usepackage{dsfont}

\newtheorem{thm}{Theorem}[section]
\newtheorem{cor}[thm]{Corollary}
\newtheorem{lem}[thm]{Lemma}
\newtheorem{defn}[thm]{Definition}
\newtheorem{prop}[thm]{Proposition}

\newtheorem{conj}[thm]{Conjecture}

\newtheorem{rem}[thm]{Remark}

\newcommand{\del}[2]{{}}

\title{Cohomology of Split Group Extensions and Characteristic Classes}

\author{Nansen Petrosyan}
\address{Department of Mathematics, Catholic University of Leuven, Kortrijk, Belgium}%
\email{Nansen.Petrosyan@kuleuven-kortrijk.be}%

\thanks{}%
\subjclass{}%
\keywords{spectral sequence, group cohomology}%

\date{\today}

\begin{document}
\begin{abstract}
There are characteristic classes that are the obstructions to the
vanishing of the differentials in the Lyndon-Hochschild-Serre
spectral sequence of a split extension of an integral lattice $L$
by a group $G$. These characteristic classes exist in the $r$-th
page of the spectral sequence provided that the differentials $d_i=0$
for all $i<r$. When $L$ decomposes into a sum of $G$-sublattices,
we show that there are defining relations between the
characteristic classes of $L$ and the characteristic classes of
its summands.
\end{abstract}

\maketitle

\section{Introduction}
Suppose $G$ is a group and $L$ is a finite rank integral $\mathbb ZG$-lattice.
Let $\Gamma = L\rtimes G$ be the semidirect product group induced by the action of $G$ on $L$.
Our objective is to analyze the split group extension $0\rightarrow L
\rightarrow \Gamma \rightarrow G\rightarrow 1$ and its associated
Lyndon-Hochschild-Serre spectral sequence $\{E_*, d_*\}$. In \cite{AGPP}
we showed that when $G$ is a cyclic group of prime
order, the spectral sequence with integral coefficients collapses at $E_2$ without
extension problems. In fact, our proof stated in \cite{AGPP} applies not only to integral
coefficients, but to all coefficient modules in
a certain category which we denote by ${\mathcal M}_F(L,G)$.  This
 category is essential in the definition of characteristic classes (see Theorem 2.2).
 It consists of all finite rank $F\Gamma$-lattices $M$ on which
$L$ acts trivially, where $F$ is a given principal ideal domain.
In particular, $F$ as a trivial $F\Gamma$-module is an object of ${\mathcal M}_F(L,G)$.

Let us now suppose that $L$ decomposes into a sum of
$\mathbb ZG$-sublattices; $L=L'\oplus L''$. Let $\{E_*^{'}, d_*^{'}\}$ and
$\{E_*^{''}, d_*^{''}\}$ be the Lyndon-Hochschild-Serre
spectral sequences associated to the respective semidirect product
groups $L^{'}\rtimes G$ and $L^{''}\rtimes G$. The
goal of the present paper is to
find necessary and sufficient  conditions
on the latter defined spectral sequences such that $\{E_*, d_*\}$
has no nonzero differentials up to a given page.

In section 2, we define characteristic classes $v_r^t(L)$, introduced by
Charlap-Vasquez  and  Sah (see \cite{CV} and \cite{sah}),
that are the obstructions to the vanishing of the differentials in
$\{E_*, d_*\}$. These classes exist and lie in the images of the
differentials $d_r^{0,t}: E_r^{0,t}(H_t(L,F))\to E_r^{r,t-r+1}(H_t(L,F))$
provided that the differentials $d_i^{*,t}$ vanish
for all $i<r$ and for all coefficient modules in ${\mathcal M}_F(L,G)$.
Our main result is the following decomposition theorem.

\begin{thm}{\normalfont Let $G$ be any group. Assume $L'$
and $L''$ are $\mathbb ZG$-lattices of finite rank and $L=L'\oplus L''$. Let $r\geq 2$ and $t\geq 0$.

\begin{enumerate}
\item  Suppose $v^{i}_k(L')=v^{j}_k(L'')=0$ for all $i,j\leq t$ and for all $k<r$.
Then $d_2^{s,m}= \dots = d_{r-1}^{s,m}=0$ for all $s\geq 0$, all $m\leq t$,
and all coefficient modules in ${\mathcal M}_F(L,G)$. Additionally,
$$v_r^{t}(L)=\sum_{i+j=t}\Big((C^{'r}_{j})_{\ast}(v_r^{i}(L'))+(-1)^i
(C^{''r}_{i})_{\ast}(v_r^{j}(L''))\Big).$$

\smallskip

\item Suppose $v^t_{k}(L)= 0$ for all $k<r$. Then
$d^{'s,t}_{k}= d^{''s,t}_{k}= 0$ for all $s\geq 0$, all $k<r$,
and all coefficient modules in ${\mathcal M}_F(L,G)$. Additionally,
$$v^t_{r}(L')=D^{'r}_*(v^t_r(L))\hspace{3mm} \mbox{and}
\hspace{3mm} v^t_{r}(L'')=D^{''r}_*(v^t_r(L)).$$
\end{enumerate}}
\end{thm}

\medskip

In section 3, we define the $FG$-equivariant homomorphisms
$C^{'r}_{j}$, $C^{''r}_{i}$, $D^{'r}$, and $D^{''r}$
that induce the  maps $(C^{'r}_{j})_{\ast}:E_2^{'s,i-r+1}\to E_2^{s,t-r+1}$,
$(C^{''r}_{i})_{\ast}:E_2^{''s,j-r+1}\to E_2^{s,t-r+1}$, $D^{'r}_*:E_2^{s,t-r+1}\to E_2^{'s,t-r+1}$,
 and $D^{''r}_*:E_2^{s,t-r+1}\to E_2^{''s,t-r+1}$, respectively.
 The homomorphisms $C^{'2}_{j}$ and $C^{''2}_{i}$ were first
considered by Charlap and Vasquez in \cite{CV}. They showed that
the sum formula holds when $r=2$. We use a different
approach to generalize this result to all pages of the spectral sequence
and also to prove a converse.

The problem of establishing the collapse of $\{E_*, d_*\}$ in general
can be a difficult one. The spectral sequence can have
nonzero differentials even when $G$ is abelian. For instance,
Totaro proved in \cite{tot} that for any prime number $p$ and $G= {C_p^2}$,
there is a semidirect product group $\Gamma= L\rtimes {{C_p^2}} $ such that
in the associated Lyndon-Hochschild-Serre spectral sequence with $\mathbb Z_p$
coefficients there always exist nonzero differentials at $E_p$ or later.

A question, first posed by Adem (see \cite{AGPP}), that is still open is whether the spectral sequence
collapses integrally at $E_2$ without extension problems when $G$ is an arbitrary finite cyclic group. In view
of our results, we can make the following

\begin{conj}{\normalfont Let $C_n$ be a cyclic group of
order $n$ and let $L$ be a finite rank integral lattice. The Lyndon-Hochschild-Serre
spectral sequence $\{E_*, d_*\}$ of any split group extension $0\to L\to \Gamma \to C_n\to 1$
 collapses at $E_2$ for all coefficient modules in ${\mathcal M}_F(L,C_n)$.}
\end{conj}

Note that part (a) of Theorem 1.1 reduces this conjecture
to the case where $L$ is an indecomposable $\mathbb ZC_n$-lattice (see Corollary 4.2).

\section{Preliminary Results}
Henceforth, let $G$ be any group and let $L$ be any finite rank
integral $\mathbb ZG$-lattice. Denote by $\Gamma$  the associated semidirect product $L\rtimes G$.
Suppose $F$ is a principal ideal domain. For each $F\Gamma$-module $M$ we have the Lyndon-Hochschild-Serre spectral sequence
$$E^{p,q}_2(M)=H^p(G,H^q(L, M))\Longrightarrow H^{p+q} (\Gamma, M).$$
\smallskip
\indent In the proof of our main result, we make use of the multiplicative structure of this spectral sequence.
Namely, any $\Gamma$-pairing
of $F\Gamma$-modules $A$, $B$, and $C$
$$\cdot :A\otimes_F B\to C,$$
\noindent determines an $ F$-pairing
$$E_r^{p,q}(A)\otimes_F E_r^{s,t}(B)\to E_r^{p+s,q+t}(C).$$
\smallskip
In addition, for the differential $d_r$ and for each $a\in E_r^{p,q}(A)$ and
$b\in E_r^{s,t}(B)$ we have the product formula
$$d^{p+s,q+t}_r(a\cdot b)= d^{p,q}_r(a)\cdot b + (-1)^{p+q}a\cdot d^{s,t}_r(b).$$
\smallskip

We denote by ${\mathcal M}_F(L,G)$ the category of
all $F\Gamma $-modules $M$, such that $M$ is a finitely generated free
$F$-module and $L$ acts trivially on $M$. In particular, $F$ with a trivial $\Gamma$-action
is a module in this category. ${\mathcal M}_F(L,G)$ has the property
that if $M\in {\mathcal M}_F(L,G)$, then $H_i(L,M)$ and
$H^i(L,M) \cong $ Hom$_F(H_i(L, F),M)$ are also in ${\mathcal M}_F(L,G)$.
Observe that the projection of $\Gamma$ onto $G$ induces an equivalence between this category
and the category of all finite rank $FG$-lattices.

Let $t \geq 0$ and $r\geq 2$. Set $M=H_t(L,F)\cong \wedge^t(L)\otimes F$. Assume
$d_k^{0,t}=0$ for all $2\leq k\leq r-1$. By applying the Universal Coefficient Theorem it follows

\begin{align*}
E_r^{0,t}(H_t(L,F)) &= E_2^{0,t}(H_t(L,F)) \\
                    &= H^0\big(G, H^t(L, H_t(L, F))\big) \\
                    &\cong H^0\big(G, \mbox{Hom}_F(H_t(L,F), H_t(L,F))\big)\hspace{6mm}\mbox{(by UCT)}\\
                    &= \mbox{Hom}_{FG}\big(H_t(L,F), H_t(L,F)\big)\\
                    &\cong \mbox{Hom}_{FG}\big(\wedge^t(L)\otimes F, \wedge^t(L)\otimes F\big).\\
  \end{align*}
\begin{defn}{\normalfont With the preceding assumptions, let $\mbox{id}^t:\wedge^t(L)\otimes F\to \wedge^t(L)\otimes F$
be the identity homomorphism. Identifying along the above isomorphisms, denote by $[f]$  the class in $E_r^{0,t}(H_t(L,F))$
corresponding to a map $f \in \mbox{Hom}_{FG}(\wedge^t(L)\otimes F, \wedge^t(L)\otimes F)$.  Then
$v_r^t(L):=d_r^{0,t}([\mbox{id}^t])\in E_r^{r,t-r+1}(H_t(L,F))$ is said to be a \textsl{characteristic class} of
the spectral sequence $\{E_r, d_r\}$.}
\end{defn}

Characteristic classes were first considered by
Charlap and Vasquez in \cite{CV}, but only in the case when $r = 2$.
Sah in \cite{sah} extended their definition  to all $r\geq 2$
by proving the following key theorem.
\medskip

\begin{thm}\label{T:2.2}{\normalfont (Sah, 1972, \cite{sah})
Let $\{E_*, d_*\}$ be the Lyndon-Hochschild-Serre
spectral sequence of the extension $0\to L\to \Gamma \to G\to 1$. Suppose there exist
integers $r\geq 2$ and $t\geq 0$ such that $d_2^{s,t}= \dots =
d_{r-1}^{s,t}=0$ for all $s\geq 0$ and for all coefficient modules
in ${\mathcal M}_F(L,G)$.

\begin{enumerate}

\smallskip

\item  There is a canonical epimorphism
\begin{center}
$\theta : E_r^{s,0}(H^t(L,M))\rightarrow E_r^{s,t}(M)$ for all $s\geq 0$
and for all $M\in {\mathcal M}_F(L,G)$.
\end{center}

\smallskip

\item  We have
$$d_r^{s,t}(x) =(-1)^s y\cdot v_r^t(L),$$
for all $x \in E_r^{s,t}(M)$, all $M\in {\mathcal M}_F(L,G)$, and
 all $y\in E_r^{s,0}(H^t(L,M))$ with $\theta(y)=x$.

%
%

\smallskip

\item  Let $\sigma : H \rightarrow G$ be a group homomorphism which
 converts $L$ into a $\mathbb ZH$-module. Assume
$d_2^{s,t}= \dots = d_{r-1}^{s,t}=0$ holds for all $s\geq 0$ and for all
objects in ${\mathcal M}_F(L,H)$. The characteristic class
$w_r^t(L)$ for the category ${\mathcal M}_F(L,H)$ is then the
image of $v_r^t(L)$ under the map induced on the spectral sequences.
\end{enumerate}}
\end{thm}

\medskip

Note that, with the assumptions of the theorem, characteristic classes $v_r^t(L)$
are obstructions to the vanishing of the differentials $d_r^{s,t}$ for
all integers $s\geq 0$.  The following corollary is an immediate
consequence of Sah's theorem.

\begin{cor}\label{LHS}{\normalfont Let $t\geq 0$. The Lyndon-Hochschild-Serre
spectral sequence $\{E_*, d_*\}$ has $d^{s,t}_r=0$ for all $s\geq 0$, all
$r<n$, and all $M\in{\mathcal M}_F(L,G)$ if and only if the edge
differentials $d_r^{0,t}:E_r^{0,t}(H_t(L,F))\to E_r^{r,t-r+1}(H_t(L,F))$
are zero for all $r<n$.}
\end{cor}

\begin{cor}\label{LHS} {\normalfont Suppose $\varphi: L\to \Gamma$ is the
natural inclusion. Given an integer $t\geq 0$, the following statements are equivalent.

\begin{enumerate}

\item  In $\{E_*, d_*\}$, the differentials $d^{s,t}_r$ vanish for all $r,s\geq 0$ and all $M\in {\mathcal M}_F(L,G)$.

\smallskip

\item $\varphi^\ast:H^t(\Gamma, M)\rightarrow
H^t(L,M)$  maps onto  $H^t(L,M)^G$ for all $M\in {\mathcal M}_F(L,G)$.

\smallskip

\item $\varphi^\ast:H^t(\Gamma, H_t(L,F))\rightarrow
H^t(L,H_t(L,F))$  maps onto the $G$-invariants $H^t(L,H_t(L,F))^G$.
\end{enumerate}}
\end{cor}
\begin{proof}Clearly (b) implies (c). Note that the map $\varphi^\ast$ is given by the composition
$$H^t(\Gamma,M)\twoheadrightarrow E_{\infty}^{0,t}(M)=E_{t+2}^{0,t}(M)
\subset \dots \subset E_2^{0,t}(M)= H^t(L,M)^G\subset H^t(L,M).$$

\noindent If $d^{s,t}_r=0$ for all $r,s\geq 0$
and all $M\in {\mathcal M}_F(L,G)$, then
$E_{\infty}^{0,t}(M)=E_{t+2}^{0,t}(M)= \dots = E_2^{0,t}(M)$ for all $M\in {\mathcal M}_F(L,G)$.
Hence, (a) implies (b). To prove that (a) follows from (c), assume $\varphi^*:H^t(\Gamma, H_t(L,F))\rightarrow
H^t(L,H_t(L,F))^G$ is onto. Then by the above composition
we have $E_{\infty}^{0,t}(H_t(L,F))=E_{t+2}^{0,t}(H_t(L,F))= \dots = E_2^{0,t}(H_t(L,F))$.
This shows that $d_r^{0,t}:E_r^{0,t}(H_t(L,F))\to E_r^{r,t-r+1}(H_t(L,F))$
is zero for all $r\geq 2$. By the previous corollary
$d^{s,t}_r=0$ for all $s\geq 0$, $r\geq 2$, and $M\in {\mathcal M}_F(L,G)$.
\end{proof}

\section{Decomposition Theorem} Suppose $L$
is a direct sum of $\mathbb ZG$-sublattices $L'$ and $L''$.
In this section we derive relations between the characteristic classes of
$L$ and the characteristic classes of $L'$ and $L''$.

For convenience, we use $\wedge^*(L)$ to denote $\wedge ^*(L)\otimes F$.
Recall that there is a standard $FG$-module decomposition
$$\wedge^n(L)\cong
\bigoplus_{i+j=n}\wedge^i(L')\otimes\wedge^j(L'').$$

\smallskip

\begin{defn}{\normalfont Given  integers $i, j\geq 0$ and $r\geq 1$, define
$FG$-equivariant homomorphisms

\begin{align*}
&C^{'r}_{j}:\mbox{Hom}_{F}(\wedge^i(L'),\wedge^{i+r-1}(L'))\rightarrow
\mbox{Hom}_{F}(\wedge^{i+j}(L),\wedge^{i+j+r-1}(L))\\
\intertext{by} &C^{'r}_{j}(f)(x\otimes y)= \left\{\begin{array}{ll} f(x)\otimes y &
\mbox{ if $x\in \wedge^i(L')$ and $ y\in \wedge^j(L'')$ }\\0 &
\mbox{ otherwise}\\
\end{array}
\right.\\
\intertext{and}
&C^{''r}_{i}:\mbox{Hom}_{F}(\wedge^j(L''),\wedge^{j+r-1}(L''))\rightarrow
\mbox{Hom}_{F}(\wedge^{i+j}(L),\wedge^{i+j+r-1}(L))\\
\intertext{by} &C^{''r}_{i}(f)(x\otimes y)= \left\{\begin{array}{ll} x\otimes f(y) &
\mbox{ if $x\in \wedge^i(L')$ and $ y\in \wedge^j(L'')$ }\\0 &
\mbox{ otherwise}.\\
\end{array}\right.
\end{align*}}

\end{defn}


Let $\Gamma'=L'\rtimes G$ and $\Gamma''=L''\rtimes G$. There are associated
split exact sequences of groups
$$0\rightarrow L'\rightarrow \Gamma'\rightarrow G \rightarrow 1 \hspace{5mm} \mbox{and}
\hspace{5mm} 0\rightarrow L''\rightarrow \Gamma'' \rightarrow G \rightarrow 1.$$

\smallskip

\noindent Suppose
$v^*_\ast(L')$ and $v^*_\ast(L'')$ are the respective
characteristic classes of the Lyndon-Hochschild-Serre spectral
sequences $\{E'_{\ast}, d^{'}_\ast\}$ and $\{E''_{\ast}, d^{''}_\ast\}$
corresponding to these extensions. Let $\iota': L'\to L$ be
the natural inclusion and let ${p}':L\to L'$ be the natural projection.
Similarly, define $\iota'': L''\to L$ and ${p}'':L\to L''$.
As a straightforward application of the definitions, we obtain the following lemma.

\smallskip

\begin{lem}{\normalfont The map
$C^{'r}_0:\mbox{Hom}_F(\wedge^i(L'),\wedge^{i+r-1}(L'))\rightarrow
\mbox{Hom}_F(\wedge^{i}(L),\wedge^{i+r-1}(L))$ is  given by the
composition ${p}'^*\circ \iota'_*$, where
\begin{align*}
&\mbox{Hom}_F(\wedge^i(L'),\wedge^{i+r-1}(L'))\buildrel{\iota'_{*}}\over{\longrightarrow}
\mbox{Hom}_F(\wedge^{i}(L'),\wedge^{i+r-1}(L))\buildrel{{p}'^*}\over{\longrightarrow}
\mbox{Hom}_F(\wedge^{i}(L),\wedge^{i+r-1}(L)),\\
\intertext{and $C^{''r}_0$ is the composition ${p}''^*\circ \iota''_*$, where}
&\mbox{Hom}_F(\wedge^j(L''),\wedge^{j+r-1}(L''))\buildrel{\iota''_{*}}\over{\longrightarrow}
\mbox{Hom}_F(\wedge^{j}(L''),\wedge^{j+r-1}(L))\buildrel{{p}''^*}\over{\longrightarrow}
\mbox{Hom}_F(\wedge^{j}(L),\wedge^{j+r-1}(L)).
\end{align*}}
\end{lem}

\medskip

\begin{proof}Suppose $f \in \mbox{Hom}_F(\wedge^i(L'),\wedge^{i+r-1}(L'))$.  For any $x\in \wedge^i(L')$
and $y\in \wedge^0(L'')$,  $(p'^* \circ \iota'_*) (f)(x\otimes y)=y
(\iota'_*(f)(x))=y(f(x)\otimes 1)=f(x)\otimes y=C^{'r}_0(f)(x\otimes
y).$ The second assertion of the lemma follows analogously.
\end{proof}

\smallskip

\begin{defn}{\normalfont Given integers $i,j\geq 0$ and $r\geq 1$, let
$\wedge^{i+r-1}p': \wedge^{i+r-1}(L)\rightarrow \wedge^{i+r-1}(L')$ and
$\wedge^{j+r-1}p'': \wedge^{j+r-1}(L)\rightarrow \wedge^{j+r-1}(L'')$ be the maps
induced by the projections $p'$ and $p''$, respectively. Define $FG$-equivariant homomorphisms

\begin{align*}
&D^{'r}:\mbox{Hom}_F(\wedge^i(L),\wedge^{i+r-1}(L))\rightarrow
\mbox{Hom}_F(\wedge^{i}(L'),\wedge^{i+r-1}(L'))\\
\intertext{by}
&D^{'r}(f)(x)=  \wedge^{i+r-1}p'(f(x\otimes 1)) \mbox{ for all } x\in \wedge^{i}(L'),\\
\intertext{and}
&D^{''r}:\mbox{Hom}_F(\wedge^j(L),\wedge^{j+r-1}(L))\rightarrow
\mbox{Hom}_F(\wedge^{j}(L''),\wedge^{j+r-1}(L''))\\
\intertext{by}
&D^{''r}(f)(y)=  \wedge^{j+r-1}p''(f(1\otimes y)) \mbox{ for all } y\in \wedge^{j}(L'').
\end{align*}}
\end{defn}

\begin{lem}{\normalfont The map $D^{'r}$ is given by the composition ${p}'_*\circ \iota'^*$, where
\begin{align*}
&\mbox{Hom}(\wedge^i(L),\wedge^{i+r-1}(L))\buildrel{\iota'^{\ast}}\over{\longrightarrow}
\mbox{Hom}(\wedge^{i}(L'),\wedge^{i+r-1}(L))\buildrel{p'_{\ast}}\over{\longrightarrow}
\mbox{Hom}(\wedge^{i}(L'),\wedge^{i+r-1}(L')),\\
\intertext{and $D^{''r}$ is the composition ${p}''_*\circ \iota''^*$, where}
&\mbox{Hom}(\wedge^j(L),\wedge^{j+r-1}(L))\buildrel{\iota''^{\ast}}\over{\longrightarrow}
\mbox{Hom}(\wedge^{j}(L''),\wedge^{j+r-1}(L))\buildrel{p''_{\ast}}\over{\longrightarrow}
\mbox{Hom}(\wedge^{j}(L''),\wedge^{j+r-1}(L'')).
\end{align*}}
\end{lem}
\begin{proof}This is an immediate consequence of the definitions of $D^{'r}$ and $D^{''r}$.
\end{proof}

\smallskip

\begin{prop} {\normalfont Let $i\geq 0$ and $r\geq 2$.
Suppose $E^{'s,t}_r(M)=E^{'s,t}_2(M)$
and $E^{s,t}_r(H_i(L,F))=E^{s,t}_2(H_i(L,F))$  when
$(s,t)=(0,i)$ and when $(s,t)=(r,i-r+1)$, for $M=H_i(L',F)$ and for
$M=H_i(L,F)$. Then $d^{0,i}_r\circ
(C^{'1}_{0})_*=(C^{'r}_{0})_*\circ d^{'0,i}_r$ and $d^{'0,i}_r\circ
D^{'1}_*=D^{'r}_*\circ d^{0,i}_r$.}
\end{prop}

\begin{proof}The first claim asserts that the following diagram commutes.

$$\begin{CD} H^0\big(G,H^i(L',H_i(L', F))\big) @ >d'^{0,i}_r>> H^r\big(G,H^{i-r+1}(L',H_i(L', F))\big)\\
@V{\operatorname (C^{'1}_{0})_*}VV @V{\operatorname (C^{'r}_{0})_*}VV @.\\
H^0\big(G,H^i(L,H_i(L, F))\big) @ >d^{0,i}_r>> H^r\big(G,H^{i-r+1}(L,H_i(L,
F))\big)
\end{CD}$$

\bigskip


\noindent By Lemma 3.2, we know that the map $(C^{'*}_{0})_*$ is
induced by the inclusion $\iota':L'\to L$ and the projection $p':L\to L'$.
Hence, this diagram is the outer square of the commutative
diagram below.

$$\begin{CD} H^0\big(G,H^i(L',H_i(L', F))\big) @ >d'^{0,i}_r>> H^r\big(G,H^{i-r+1}(L',H_i(L', F))\big)\\
@V{\operatorname \iota'_*}VV @V{\operatorname \iota'_*}VV @.\\
H^0\big(G,H^i(L',H_i(L, F))\big) @ >d'^{0,i}_r>> H^r\big(G,H^{i-r+1}(L',H_i(L,
F))\big)\\
@V{p'^*}VV @V{p'^*}VV @.\\
H^0\big(G,H^i(L,H_i(L, F))\big) @ >d^{0,i}_r>> H^r\big(G,H^{i-r+1}(L,H_i(L,
F))\big)
\end{CD}$$

\bigskip

\noindent The second claim follows by a similar argument using Lemma 3.4.
\end{proof}

\smallskip

As previously noted, the category ${\mathcal M}_F(L,G)$ can be identified
with the category of all finite rank $FG$-lattices. In view of this, we
will not distinguish
between  ${\mathcal M}_F(L',G)$, ${\mathcal M}_F(L'',G)$, and ${\mathcal M}_F(L,G)$.

\begin{lem}{\normalfont There are natural morphisms of spectral sequences  $\pi'^*:E_r^{'}\to E_r$,
 $\phi'^*:E_r\to E_r^{'}$ and $\pi''^*:E_r^{''}\to E_r$,
 $\phi''^*:E_r\to E_r^{''}$, such that  for all $r\geq 2$, all $s\geq 0$, and
 all coefficient modules in ${\mathcal M}_F(L,G)$,
$$d_r^{'s,t}= \phi'^* \circ d_r^{s,t}\circ \pi'^* \hspace{3mm} \mbox{and}
\hspace{3mm} d_r^{''s,t} = \phi''^* \circ d_r^{s,t}\circ \pi''^*.$$}
\end{lem}

\begin{proof}We observe that the map $\iota'$ induces a natural inclusion  $\phi' : \Gamma'\to \Gamma$
and the map $p'$ induces a natural projection $\pi':\Gamma\to \Gamma'$ such that the composition
$\pi'\circ \phi'$ is the identity map on $\Gamma'$. Let $\varphi:L\to \Gamma$ and $\varphi':L'\to \Gamma'$ be the
canonical inclusions.
 It follows that $\phi'\circ \varphi'= \varphi\circ \iota'$ and $\pi'\circ \varphi = \varphi' \circ p'$.
Hence, the homomorphisms $\phi'$ and $\pi'$ give rise to spectral
sequence morphisms $\phi'^*:E_*\to E_*^{'}$ and $\pi'^*:E_*^{'}\to E_*$, such that $\phi'^*\circ \pi'^*$ is
the identity morphism on $\{E_*^{'}, d_*^{'}\}$. Then, $d_r^{'s,t}= \phi'^* \circ \pi'^* \circ d_r^{'s,t}
= \phi'^* \circ d_r^{s,t}\circ \pi'^*$. Analogously, $\iota''$ and $p''$ induce an inclusion  $\phi'' : \Gamma''\to \Gamma$
and a projection $\pi'':\Gamma\to \Gamma''$, respectively, such that the composition
$\pi''\circ \phi''$ is the identity map on $\Gamma''$ and $d_r^{''s,t}= \phi''^* \circ d_r^{s,t}\circ \pi''^*$.
\end{proof}

\medskip

We are now ready to prove our main result.

\begin{thm}{\normalfont Recall that $G$ is an arbitrary group, and that $L$ is a $\mathbb ZG$-lattice
of finite rank that decomposes into a direct sum of the $\mathbb ZG$-lattices $L'$ and $L''$. Let $r\geq 2$ and $t\geq 0$.

\begin{enumerate}
\item  Suppose $v^{i}_k(L')=v^{j}_k(L'')=0$ for all $i,j\leq t$ and for all $k<r$.
Then $d_2^{s,m}= \dots = d_{r-1}^{s,m}=0$ for all $s\geq 0$, all $m\leq t$,
and all coefficient modules in ${\mathcal M}_F(L,G)$. Additionally,
$$v_r^{t}(L)=\sum_{i+j=t}\Big((C^{'r}_{j})_{\ast}(v_r^{i}(L'))+(-1)^i
(C^{''r}_{i})_{\ast}(v_r^{j}(L''))\Big).$$

\smallskip

\item Suppose $v^t_{k}(L)= 0$ for all $k<r$. Then
$d^{'s,t}_{k}= d^{''s,t}_{k}= 0$ for all $s\geq 0$, all $k<r$,
and all coefficient modules in ${\mathcal M}_F(L,G)$. Additionally,
$$v^t_{r}(L')=D^{'r}_*(v^t_r(L))\hspace{3mm} \mbox{and}
\hspace{3mm} v^t_{r}(L'')=D^{''r}_*(v^t_r(L)).$$
\end{enumerate}}
\end{thm}

\bigskip

\begin{proof}

To show part (a), we will use induction on $k\geq 2$ to prove that $v_k^{m}(L)$ satisfies the sum
formula for all $m\leq t$ and all $k\leq r$. For each $k< r$,
 since $v^{i}_{k}(L')=v^{j}_k(L'')=0$ for all $i,j\leq t$,
this will imply  $v_k^{m}(L)=0$ and hence, by Theorem 2.2(b), $d_k^{s,m}=0$ for all $s\geq 0$,
all $m\leq t$, and all $M\in{\mathcal M}_F(L,G)$.

Suppose $d_2^{s,m}= \cdots = d_{k-1}^{s,m}=0$ for all $s\geq 0$,
all $m\leq t$, and all $M\in {\mathcal M}_F(L,G)$. Recall that
 $v^m_k(L) =
d_k^{0,m}([\mbox{id}^m])$, where $\mbox{id}^m\in \mbox{Hom}_{FG}(\wedge^m(L),\wedge^m(L))$
is the identity map. Similarly, $v^{i}_k(L') =
d_k^{'0,i}([\mbox{id}'^{i}])$ and $v^{j}_k(L'') = d_k^{''0,j}([\mbox{id}''^{j}])$.
Consider the decomposition  $\mbox{id}^{m}=\displaystyle{\sum_{i+j=m}}\mbox{id}_{ij}$, where
$\mbox{id}_{ij}:\wedge^{i+j}(L)\rightarrow \wedge^{i+j}(L)$ is the $F$-linear map given by
$$\mbox{id}_{ij}(x\otimes y)= \left\{\begin{array}{ll} x\otimes y & \mbox{ if
$x\in \wedge^i(L')$ and $ y\in \wedge^j(L'')$ }\\0 & \mbox{
otherwise.}\\\end{array}\right.$$

\bigskip

\noindent Given $i,j\geq 0$, let $x\in \wedge^i(L')$ and $y\in
\wedge^j(L'')$. Then, $\mbox{id}_{ij}(x\otimes y)= x\otimes y = (x\otimes
1)\wedge (1\otimes y) =C'^1_0(\mbox{id}'^i)(x\otimes 1)\wedge
C''^1_0(\mbox{id}''^j)(1\otimes y).$ This implies $[\mbox{id}_{ij}]=
(C'^1_{0})_*([\mbox{id}'^i])\cdot (C''^1_{0})_*([\mbox{id}''^j])\in H^0(G,
H^{i+j}(L,H_{i+j}(L,F)))$, and thus

$$[\mbox{id}^{m}]=\sum_{i+j=m}(C'^1_{0})_*([\mbox{id}'^i])\cdot
(C''^1_{0})_*([\mbox{id}''^j]).$$

\medskip

\noindent By applying the product formula for the differentials, we
compute

\begin{align*}
v^{m}_k(L) &=d_k^{0,m}\Big(\sum_{i+j=m}(C'^1_{0})_*([\mbox{id}'^i])\cdot
                       (C''^1_{0})_*([\mbox{id}''^j])\Big) \\
                      &=\sum_{i+j=m}d_k^{0,m}\Big((C'^1_{0})_*([\mbox{id}'^i])\cdot
                       (C''^1_{0})_*([\mbox{id}''^j])\Big)\\
                      &=\sum_{i+j=m}\Big(d_k^{0,i}\big((C'^1_{0})_*([\mbox{id}'^i])\big)\cdot
                       (C''^1_{0})_*([\mbox{id}''^j]) + (-1)^i (C'^1_{0})_*([\mbox{id}'^i])\cdot
                       d_k^{0,j}\big((C''^1_{0})_*([\mbox{id}''^j])\big)\Big)\\
                      &\displaystyle{\buildrel{3.5}\over {=}}\sum_{i+j=m}\Big((C'^k_{0})_*\big(d_k^{'0,i}([\mbox{id}'^{i}])\big) \cdot
                       (C''^1_{0})_*([\mbox{id}''^j]) + (-1)^i (C'^1_{0})_*([\mbox{id}'^{i}]) \cdot
                        (C''^k_{0})_*\big(d_k^{''0,j}([\mbox{id}''^{j}])\big)\Big)\\
                      &=\sum_{i+j=m}\Big((C'^k_{0})_*(v_k^{i}(L'))\cdot
                       (C''^1_{0})_*([\mbox{id}''^j]) + (-1)^i (C'^1_{0})_*([\mbox{id}'^{i}])\cdot
                        (C''^k_{0})_*(v_k^{j}(L''))\Big)\\
                      &=\sum_{i+j=m}\Big((C'^k_{j})_*(v_k^{i}(L'))+(-1)^i
                        (C''^k_{i})_*(v_k^{j}(L''))\Big).\\
\end{align*}
\noindent The last equality follows from the definitions of
$C^{'\ast}_{\ast}$ and $C^{'' \ast}_{\ast}$ and the fact that
$\mbox{id}'^{i}$ and $\mbox{id}''^{j}$ are identity maps.

To prove (b), we will use induction on $k\geq 2$ to show that
$v^t_{k}(L')=D^{'r}_*(v^t_k(L))$ and $v^t_{k}(L'')=D^{''r}_*(v^t_k(L))$ for
all $k\leq r$. By Theorem 2.2, for each $k< r$ and for all $s\geq 0$ this
will imply that $d_k^{'s,t}=d_k^{''s,t}=0$ for all coefficient modules
in ${\mathcal M}_F(L,G)$. Note that this assertion also follows from Lemma 3.6.

For the identity map
$\mbox{id}'^t:\wedge^t(L')\to \wedge^t(L')$, we have $\mbox{id}'^t(x)=x=\wedge^tp'(x\otimes 1)=\wedge^tp'(\mbox{id}^t(x\otimes 1))=
D^{'1}(\mbox{id}^t)(x)$ for all  $x\in \wedge^t(L')$. This implies $[\mbox{id}'^{t}]= D^{'1}_*([\mbox{id}^t])\in E^{'0,t}_k(H_t(L',F))$.
Thus, it follows
$$
               v_k^t(L') =d_k^{'0,t}([\mbox{id}'^{t}])
                          =d_k^{'0,t}\circ D^{'1}_*([\mbox{id}^{t}])
                         \displaystyle{ \buildrel{3.5}\over { = } } D^{'k}_*\circ d_k^{0,t}([\mbox{id}^{t}])
                           =D^{'k}_*(v^t_k(L)).
$$
\noindent By an analogous argument,  $v^t_{k}(L'')=D^{''k}_*(v^t_k(L))$.
\end{proof}

\smallskip

\begin{rem}{\normalfont Using the same assumptions as in Theorem 3.7(a), the sum formula can be simplified.
Since $E_r^{'r, i-r+1}=E_r^{''r, j-r+1}=0$, $v^{i}_r(L')= v^{j}_r(L'')= 0$ holds a priori for all $i,j< r-1$.
It is an easy exercise to check that $d_r^{'0,r-1}=0$ and $d_r^{''0,r-1}=0$
for all $M\in {\mathcal M}_F(L,G)$, since they are differentials with target in the $0$th row (see Prop.1, \cite{sah}).
Therefore, $v^{r-1}_r(L')= v^{r-1}_r(L'')= 0$ and we have
\begin{align*}
v^t_r(L)&=\sum_{i+j=t}\Big((C^{'r}_{j})_*(v_r^{i}(L'))+(-1)^i
(C^{''r}_{i})_*(v_r^{j}(L''))\Big)\\
        &=\sum_{i+j=t-r}(C^{'r}_{j})_*(v_r^{i+r}(L')) + \sum_{i+j=t-r}(-1)^i
          (C^{''r}_{i})_*(v_r^{j+r}(L''))\\
        &=\sum_{i+j=t-r}\Big((C^{'r}_{j})_*(v_r^{i+r}(L'))+(-1)^i
          (C^{''r}_{i})_*(v_r^{j+r}(L''))\Big).\\
\end{align*}}
\end{rem}

\section{Some Corollaries}

Theorem 3.7 has particularly interesting applications if characteristic classes of the
$\mathbb ZG$-lattices $L'$, $L''$, and $L$ are viewed as obstructions to the vanishing
of the differentials in the associated Lyndon-Hochschild-Serre spectral sequences. We use
the same notation as before.

\begin{cor}\label{LHS}{\normalfont Let $r\geq 2$. If $v^{i}_k(L')=v^{j}_k(L'')=0$ for all
$i,j \in [k,r]$ and for all $k< r$, then $d_2^{s,m}= \dots =
d_{r-1}^{s,m}=0$ for all $s\geq 0$, all $m\leq r$, and all $M\in {\mathcal M}_{F}(L,G)$.
Moreover,
$$v_r^{r}(L)=(C^{'r}_{0})_*(v_r^{r}(L'))+ (C^{''r}_{0})_*(v_r^{r}(L'')).$$}
\end{cor}

\begin{proof}This is a direct consequence of Theorem 3.7 and the preceding remark.

\end{proof}

\smallskip

A consequence of Lemma 3.6 is the fact that when the differentials in the
Lyndon-Hochschild-Serre spectral sequence
$\{E_*, d_*\}$ are all zero, the same is true for the spectral sequences corresponding to the
$\mathbb ZG$-sublattices $L'$ and $L''$. The next corollary gives us a converse.

\begin{cor}\label{LHS}{\normalfont  Suppose
$v^t_{p}(L')=v^t_{q}(L'')=0$ for all $p \leq \mbox{dim} (L')$, all $ q \leq \mbox{dim} (L'')$, and all
$t\geq 0$. Then $d_k^{s,t}=0$ for all $M\in {\mathcal M}_{F}(L,G)$, all $s,t\geq 0$, and all $k\geq 2$.
Moreover, if in addition $\{E_*, d_*\}$ has no extension problems, then for every $n\geq 0$ and
for all $M\in {\mathcal M}_{F}(L,G)$ we have $$H^n(\Gamma, M)=\bigoplus_{i+j=n}H^i(G,H^j(L, M)).$$}
\end{cor}

\begin{proof} Since $v^t_{p}(L')=d_p^{'0,t}([\mbox{id}'^t])\in E_p^{p,t-p+1}(H_t(L',F))$, this class lies
in the image of the map
$$d_p^{'0,t}:H^0(G, H^t(L',H_t(L',F)))\to H^p(G, H^{t-p+1}(L', H_{t}(L',F))).$$
Note that $d_p^{'0,t}=0$ when $t>\mbox{dim}(L')$ or $p>t+1$. If $p=t+1$, then $d_{t+1}^{'0,t}=0$,
since it is a differential with a target in the $0$th row
(see Prop.1, \cite{sah}). Therefore, if $v^*_{p}(L')=0$ for all
$p\leq \mbox{dim}(L')$, then all characteristic classes of the spectral sequence
 $\{E'_*, d'_*\}$ are zero. A similar argument shows that
all characteristic classes of
 $\{E''_*, d''_*\}$ are zero when $v^*_{q}(L'')=0$ for all
$q\leq \mbox{dim}(L'')$.
\end{proof}

\smallskip

\begin{cor}\label{LHS}{\normalfont Let $t\geq 0$. Set $\Gamma^{'}=L'\rtimes G$ and
 $\Gamma^{''}=L''\rtimes G$. Let $\varphi':L'\to \Gamma^{'}$, $\varphi'':L''\to \Gamma^{''}$, and $\varphi:L\to \Gamma$
 be the natural inclusions.

\smallskip

\begin{enumerate}
\item If
$\varphi'^\ast:H^m(\Gamma^{'}, H_m(L',F))\rightarrow
H^m(L',H_m(L',F))^G$ and $\varphi''^\ast:H^m(\Gamma^{''},
H_m(L'',F))\rightarrow H^m(L'',H_m(L'',F))^G$ are surjective for all $m\leq t$, then
$\varphi^\ast:H^m(\Gamma, M) \rightarrow H^m(L,M)^G$ is surjective
for all $m\leq t$ and for all $M\in {\mathcal
M}_F(L,G)$.

\smallskip

\item If $\varphi^\ast:H^t(\Gamma, M) \rightarrow H^t(L,M)^G$ is surjective, then
$\varphi'^\ast:H^t(\Gamma^{'}, M)\rightarrow
H^t(L',M)^G$ and $\varphi''^\ast:H^t(\Gamma^{''},
M)\rightarrow H^t(L'',M)^G$ are surjective for all $M\in
{\mathcal M}_F(L,G)$.
\end{enumerate}}
\end{cor}

\begin{proof} To prove (a), we observe that Corollary 2.4 implies
$v^i_r(L')=v^j_r(L'')=0$ for all $r\geq 0$ and for all $i,j\leq t$. Then, by Theorem 3.7,
 $d_r^{s,m}=0$ for all $r,s\geq 0$, all $m\leq t$,
and all coefficient modules in ${\mathcal M}_F(L,G)$. Applying again Corollary 2.4 finishes
the proof.

For part (b), let $\iota'^*:H^t(L,M)^G\to H^t(L',M)^G$, $\phi'^*:H^t(\Gamma, M)\to H^t(\Gamma', M)$, and
$p'^*:H^t(L',M)^G\to H^t(L,M)^G$ be the induced maps of the inclusions $\iota':L'\to L$,
 $\phi':\Gamma'\to \Gamma$, and the projection $p':L\to L'$, respectively. Then, we have
$\varphi'^*\circ \phi'^*= \iota'^*\circ \varphi^*$. Since $\iota'^*\circ p'^*$ is the identity map on
$H^t(L',M)^G$, $\iota'^*$ is surjective and hence $\iota'^*\circ \varphi^*$ is  surjective.
The previous equality shows that $\varphi'^*$ is also surjective. Similarly, it follows that $\varphi''^*$ is surjective.
\end{proof}

\smallskip

Given an arbitrary finite rank integral $\mathbb ZG$-lattice $L$,
Lieberman's result (Thm.4, \cite{sah}) states that in the associated
Lyndon-Hochschild-Serre spectral sequence $\{E_*, d_*\}$, for any $s,t\geq 0$ and $r\geq 2$, the image of
the differential $d_r^{s,t}$
is a torsion group annihilated by the integers $m^{t-r+1}(m^{r-1}-1)$ for all $m\in \mathbb Z$.
Using this fact, it was proved in \cite{sah}  that if $F$ is a field of nonzero characteristic $p$, then
$ d_r^{s,t}=0$ for all $s,t\geq 0$ and all $r<p$. The next corollary follows from
combining this result with the sum formula of Remark 3.8.

\begin{cor}\label{LHS}{\normalfont Let  $F$ be a field of nonzero characteristic $p$.
Assume $L'$ and $L''$ are $\mathbb ZG$-lattices of finite rank and $L=L'\oplus L''$.

\smallskip

\begin{enumerate}
\item $ d_r^{s,t}=0$ for all $s,t\geq 0$, all $r<p$, and all $M\in {\mathcal M}_F(L,G)$.

\medskip

\item $ d_p^{s,t}(x) = (-1)^s y\cdot
\displaystyle{\sum_{i+j=t-p}}\Big((C^{'p}_{j})_*(v_p^{i+p}(L'))+(-1)^i
(C^{''p}_{i})_*(v_p^{j+p}(L''))\Big)$

\smallskip

\noindent for all $M\in {\mathcal M}_F(L,G)$, for all $x \in E_p^{s,t}(M)$, and
for all $y\in E_p^{s,0}(H^t(L,M))$ such that  $\theta(y)=x$.
\end{enumerate}}
\end{cor}\qed
$${\bf Acknowledgments }$$

I am thankful to Alejandro Adem for suggesting to me to study the cohomology of semidirect product
groups and for always giving  good advice. I also thank Jim Davis and Karel Dekimpe for
their many conversations. Lastly, I would like to thank the referee for many helpful comments
 and suggestions which led to a better exposition of the results in the article.


\begin{thebibliography}{12345}

\bibitem{AGPP} Adem, A., Ge, J., Pan, J., Petrosyan, N. \emph{Compatible Actions and
Cohomology of Crystallographic Groups}, J. of Algebra \textbf{320} (2008), 341--353.




\bibitem{CV} Charlap, L. Vasquez, A., \emph {Characteristic
classes for modules over groups I}, Trans. Amer. Math. Soc.
\textbf{137} (1969), 533--549.

\bibitem{sah} Sah, C.-H., \emph {Cohomology of split group extensions}, J. of Algebra
\textbf{29} (1974), 255--302.

\bibitem{tot} Totaro, B., \emph{Cohomology of Semidirect Product Groups},
J. of Algebra \textbf{182} (1996), 469--475.

\end{thebibliography}
\end{document}